\documentclass[12pt,reqno]{article}
\usepackage{amsmath, amsthm, graphicx,amsfonts, amssymb,color}
\usepackage{bm}
\usepackage{hyperref}

\usepackage{amsmath}
\numberwithin{equation}{section}
\def\d{\mathrm{d}}\def\e{\mathrm{e}}
\newtheorem{theorem}{Theorem}[section]
\newtheorem{lemma}[theorem]{Lemma}
\newtheorem{corollary}[theorem]{Corollary}

\newtheorem{remark}[theorem]{Remark}
\newtheorem{example}[theorem]{Example}

\def\<{\langle}\def\>{\rangle}
\def\prf{\noindent{\bf Proof.~~}}
\def\deprf{\hfill$\Box$\medskip}

\usepackage{mathrsfs}
\newcommand{\scr}[1]{\mathscr #1}

\def\d{\mathrm{d}}

\def\bg{\begin}

\def\de{\end{equation}}

\def\edar{\end{eqnarray}}

\def\l{\left}\def\r{\right}

\def\[{\l[} \def\]{\r]}
\def\({\l(} \def\){\r)}

\def\beqlb{\begin{eqnarray}}\def\eeqlb{\end{eqnarray}}
\def\beqnn{\begin{eqnarray*}}\def\eeqnn{\end{eqnarray*}}
\def\d{{\mbox{\rm d}}}

\title{{\bf Ergodic convergence rates for time-changed symmetric L\'{e}vy processes in dimension one} }
\author{
{\bf Tao Wang\footnote{
	E-mail: wang\_tao@mail.bnu.edu.cn. }}\\
\footnotesize{School of Mathematical Sciences, Beijing Normal University, }\\
\footnotesize{Laboratory of Mathematics and Complex Systems, Ministry of Education}\\
\footnotesize{Beijing 100875, China}\\
}
\date{ }

\date{}

\begin{document}

\maketitle


\bg{abstract}
We obtain the lower bounds for ergodic convergence rates, including spectral gaps and convergence rates in strong ergodicity, for time-changed  symmetric L\'{e}vy processes, by using harmonic function and reversible measure.
As direct applications, explicit sufficient conditions for exponential and strong ergodicity are given. Some examples are also presented. 
\end{abstract}




\section{Main results and examples}
Ergodicity for  L\'{e}vy-type processes is an important topic in  study of Markov processes.  In general,  criteria are obtained by using Lyapunov functions (cf.  \cite{wj08} for general one-dimensional L\'{e}vy-type operators, \cite{wj13}  for L\'{e}vy-driven SDEs, and \cite{CW14}  for  time-changed symmetric stable processes), or coupling methods (see \cite{wj19}  for L\'{e}vy-driven SDEs). 

Recently, \cite{W21} obtains the criteria for strong and exponential ergodicity of one-dimensional time-changed symmetric stable processes;  the lower bounds  for ergodic convergence rates, including spectral gaps and convergence rates in strong ergodicity are also estimated in \cite{W21}. Different from Lyapunov criteria and coupling methods, the main idea in \cite{W21} is to estimate the Green function for $\mathbb{R}\setminus\{0\}$ by assuming that the process is pointwise-recurrent.

While \cite{W21} deals with the classical  $\alpha$-stable processes for $\alpha\in (1,2)$, this will exclude some significant pointwise-recurrent L\'{e}vy processes such as  the  diffusion operator with stable jump: $L=a(x)(c_1\Delta+c_2\Delta^{\alpha/2}),$
where $\alpha\in(1,2)$, $a$ is a positive function such that $a^{-1}$ is Lebesgue integrable, $c_1$ and $c_2$ are two constants, $\Delta$ is the Laplacian operator and $\Delta^{\alpha/2}$ is the fractional Laplacian. 

The aim of this paper is to study  the ergodic convergence rates for general pointwise-recurrent time-changed symmetric L\'{e}vy processes. To this end, we first recall some basic definitions. 

Let $X=\left(X_{t}\right)_{t \geqslant 0}$ be a one-dimensional  symmetric  L\'{e}vy process with L\'{e}vy measure $\nu$.  
The characteristic exponent of $X$ is defined as
\begin{equation}\label{symbol}
	\psi(\xi):=-\log\mathbb{E}_0\e^{i\xi X_1}=\int_{\mathbb{R}}(1-\cos \xi x) \nu(\d x)+\sigma^{2} \xi^{2}, \quad \xi,\sigma \in \mathbb{R}.
\end{equation}
The  generator $\mathcal{A}$ is a one-dimensional L\'{e}vy operator  given by
\begin{equation}\label{generator}
	\mathcal{A} u(x):=-\int \psi(\xi) \widehat{u}(\xi) e^{i x \xi} \d \xi, \quad u \in C_{0}^{\infty}(\mathbb{R}),
\end{equation}
where $C_{0}^{\infty}(\mathbb{R})$  is the space of compactly supported smooth functions in $\mathbb{R}$ and
$\widehat{u}(\xi)$ is the Fourier transform of $u$.
The corresponding regular Dirichlet form $(\scr{E},\mathcal{F})$ is given by
\begin{equation}\label{Diri form}
	\begin{split}
		\mathscr{E}(u, w)
		&=  \int_{\mathbb{R}}  \frac{1}{2} \sigma^2u'(x) w'(x) \d x +\iint_{\mathbb{R} \times \mathbb{R} \backslash \text { diag }}(u(x+h)-u(x))(w(x+h)-w(x)) \nu(\d h) \d x,\\
	\end{split}
\end{equation}
for  $ u, w \in \mathcal{F}:=\{f\in L^2(d x):\scr{E}({f,f})<\infty\}$
(see \cite[Example 3.13]{bsw13} for more details).

Let   $p_t(y-x):=p_t(x,y)$ be the transition density of this process. 
Define
\begin{equation}\label{harmonic}
	H(x)=\int_{0}^\infty (p_s(0)-p_s(x))\d s=\frac{1}{\pi} \int_{0}^{\infty}(1-\cos x s) \frac{1}{\psi(s)} \d s.
\end{equation}
Obviously, $H$ is an even function: $H(x)=H(-x)$.
Let $P_t^{0}$ be the semigroup of the process killed upon hitting the origin, i.e.,
\begin{equation}\label{killed semigroup}
	P_t^{0}(x, A)=\mathbb{P}_x[X_t\in A,t<\tau_0],\quad \forall x\in \mathbb{R}, \ \forall A\in\mathscr{B}(\mathbb{R}),
\end{equation}
where $\tau_0:=\inf\{t>0:X_t=0\}$.
By \cite[Theorem 1.1]{ya10}, $H(x)$ is the harmonic function for $P_t^{0}$, i.e., for any $x\neq 0$, $P_t^{0}H(x)=H(x)$.

Now we consider the time-changed symmetric L\'{e}vy processes.   
Let $a$ be a positive and locally bounded measurable function on $\mathbb{R}$ with $a(x)^{-1}$  integrable. Define $\mu(\d x):=a(x)^{-1} \d x$, 
$A_{t}:=\int_{0}^{t} 1 / a(X_{s}) \d s$, and 
the
time-changed L\'{e}vy process ${Y}_{t}=X_{\tau_{t}}$, where
\begin{equation}\label{time change}
	\tau_{t}=\inf \left\{s>0: A_{s}>t\right\}.
\end{equation}
Then, the generator of time-changed process $Y$ is just the operator $L=a\mathcal{A}$ (for more details about this paragraph, see Section \ref{prem}).



Let $\mu(f):=\int_{\mathbb{R}}f\d \mu$ and $\|f\|_{L^2(\mu)}:=\sqrt{\mu(f^2)}$. We say that  $Y$  is {\bf ($L^2$-)exponentially ergodic}, if there exist non-negative function $C(x)<\infty$ and $\lambda_1>0$, such that 
\begin{equation}\label{def L2exp}
	\|P_tf-\mu(f)\|_{L^2(\mu)}\leqslant \e^{-\lambda_{1} t}\|f-\mu(f)\|_{L^2(\mu)}.
\end{equation}
The optimal convergence rate $\lambda_1$ in  \eqref{def L2exp} (i.e., the  $L^2$-spectral gap) is defined by 
\begin{equation}\label{L2-spectral}
	\lambda_{1}=\inf \{\mathscr{E}(f, f): f \in \mathscr{F}, \mu(f^{2})=1, \mu(f)=0\}.
\end{equation}

Let $\|\eta\|_{\rm{Var}}:=\sup_{|f|\leqslant 1}|\eta(f)|$ be the total variation of a signed measure $\eta$. We say that	$Y$ is {\bf strongly ergodic}, if there exist  constants $C<\infty$ and $\kappa>0$, such that
$$\sup_{x\in\mathbb{R}}\|P_t(x,\cdot)-\mu\|_{\rm{Var}}\leqslant C\mathrm{e}^{-\kappa t}.$$
The optimal convergence rate 
(see \cite{myh06} for more details)
$$\kappa= -\lim\limits_{t\rightarrow \infty}\frac{1}{t}\log \sup\limits_{x\in E}\|P_t(x,\cdot)-\mu\|_{\mathrm{Var}}=-\lim\limits_{t\rightarrow \infty}\frac{1}{t}\log \|P_t-\mu\|_{\infty\rightarrow\infty}.$$

The following theorem gives the explicit sufficient conditions for exponential ergodicity and strong ergodicity. The explicit lower bounds for ergodic convergence rates are also obtained.

\begin{theorem}\label{main thm}  Assume that $\mu(\mathbb{R})<\infty$. Consider the following conditions: 
	\begin{description}
		\item[(A1)] $\int_{0}^{\infty} ({q+\psi(x)})^{-1} \d x<\infty,\  \text{for any}\ q>0;$	
		\item[(A2)] $\int_{0}^{1} ({\psi(x)})^{-1} \d x=\infty;$	
		\item[(A3)] 
		$\psi(t)/t\rightarrow\infty$, as $t\rightarrow\infty$.
	\end{description}
	(1) If {\bf(A1)}--{\bf(A3)} hold, and
	\begin{equation}\label{delta}
		\delta:=\sup_{x} H(x)\mu\bigl((-|x|,|x|)^c\bigl)<\infty,	\end{equation}
	then  ${Y}$ is exponentially ergodic and the $L^2$-spectral gap
	$$\lambda_1\geqslant 
	\frac{1}{8\delta}.$$

	(2)	If {\bf(A1)}--{\bf(A2)} hold, and 	$$I:=\int_{\mathbb{R}}a(x)^{-1}H(|x|)d x<\infty,$$
	then   ${Y}$ is strongly  ergodic and 
	$$\kappa\geqslant\frac{1}{2I}>0.$$ 
\end{theorem}

\begin{remark}
	
	(1) For symmetric L\'{e}vy process,  condition {\bf(A1)} means that the process is not  compound Poisson and the origin is regular for itself (see \cite[Section 2]{GR17}).  Condition {\bf(A2)} means that the process is recurrent (see \cite[Section 3.2]{ya10}).  If {\bf(A1)} and {\bf(A2)} hold, then $X$ is pointwise recurrent.
	

	(2) {\bf(A3)} indicates that $H$ is  differentiable (see the proof of Theorem \ref{main thm}). 
	
	(3) Note that a time change does not change the  recurrence (cf.   \cite[Corollary 4.3.7]{OY13}). 
	Therefore, under the conditions {\bf(A1)--(A2)}, $Y$ is also pointwise recurrent, so    
	$Y$ is Lebesgue irreducible (see \cite[Page 42]{DZ96} for the definition). 
	Thus, by \cite[Theorem 4.1.1 and Theorem 4.2.1]{DZ96},   if $\mu(\mathbb{R})<\infty$, then
	$Y$ is {\it ergodic}.
	
	(4) Let $\psi(x)=|x|^{\alpha}$, $\alpha\in (1,2)$. Then  $Y$ is a time-changed symmetric $\alpha$-stable process. 
		It is well known that the process is pointwise recurrent; by \cite[Example 1.1]{ya10},  the harmonic function  $H(x)=\omega_\alpha|x|^{\alpha-1}/2$, where $\omega_\alpha=-(\cos(\pi\alpha/2)\Gamma(\alpha))^{-1}>0$. Then we have 
		$Y$ is exponentially ergodic
		if 
		\begin{equation}\label{stable-exp}
			\delta_1:=\sup_{x} |x|^{\alpha-1}\mu((-|x|,|x|)^c)<\infty,
		\end{equation}
		and
		$$\lambda_1	\geqslant
		\frac{1}{4\omega_{\alpha}\delta_1},$$
		$Y$ is strongly ergodic
		if 
		\begin{equation}\label{stable-str}
			I_1=\int_{\mathbb{R}}\sigma(x)^{-\alpha}|x|^{\alpha-1} \d x<\infty,
		\end{equation}
		and
		$$\kappa\geqslant\frac{1}{\omega_{\alpha}I_1}>0.$$ 
This case is introduced in \cite{W21}. In fact, the  conditions \eqref{stable-exp} and \eqref{stable-str} are sufficient and necessary (see \cite{W21} for more details).
\end{remark}

Next, we discuss a class of extended $\alpha$-stable processes, which is introduced in \cite{GR17}. In general, the cases
 mean that there exists $\alpha>0$ such that  $\psi(\theta)/\theta^{\alpha}$ is comparable to a non-decreasing function on $(0,\infty)$.  

We say that $\psi$ satisfies the global weak lower scaling condition, if there exist $\delta>0$  and $\beta \in(0,1]$, such that for $\lambda \geqslant 1$ and $\theta>0,$
$$
\psi(\lambda \theta) \geqslant \beta \lambda^{\delta} \psi(\theta).
$$
In short, we write $\psi\in$WLSC$(\delta,\beta)$ (see \cite{GR17} for more details). Applying this condition, we have the following result which is a direct corollary by using Theorem \ref{main thm} and \cite[Lemma 2.14]{GR17}:

\begin{corollary}\label{wlsc}
	Let $ \psi^{*}(x)=\sup _{|u| \leqslant x} \psi(u)$, $x \geqslant 0$. Assume that  there exists a constant $c>0$, such that $\psi\geqslant c\psi^*$, and $\psi\in$WLSC$(\delta,\beta)$ ($\delta>1$). If 
	$$\sup_{x} \frac{\mu\bigl((-|x|,|x|)^c\bigl)}{|x|\psi(1/x)}<\infty,$$
	then  ${Y}$ is exponentially ergodic and the spectral gap
	$$\lambda_{1}\geqslant\frac{\pi(\delta-1)\beta^2}{10}\inf_x \frac{|x|\psi^*(1/|x|)}{\mu\bigl((-|x|,|x|)^c\bigl)},$$
	
	If 	$$\int_{\mathbb{R}}\frac{1}{|x|a(x)\psi(1/x)}\d x<\infty,$$
	then   ${Y}$ is strongly  ergodic and 
	$$\kappa\geqslant\frac{\pi\beta^2(\delta-1)}{20\int_{\mathbb{R}}\bigl(|x|a(x)\psi^*(1/|x|)\bigl)^{-1}\d x}.$$
\end{corollary}



Now we  return to the diffusion operator with stable jump.

\begin{example}\label{alpha+2}
	Let $\psi(x)=c_1x^2+c_2|x|^{\alpha}$, where $c_1$, $c_2$ are two constants and $\alpha\in(1,2)$. Then $L=a(x)(c_1\Delta+c_2\Delta^{\alpha/2})$. Denote by $Y$ the corresponding process with generator $L$. By \eqref{harmonic}, the harmonic function for $c_1\Delta+c_2\Delta^{\alpha/2}$ is $$H(x)=\frac{1}{\pi} \int_{0}^{\infty} \frac{1-\cos x s}{c_1s^2+c_2s^\alpha} \d s.$$ 
	Obviously, {\bf(A1)}--{\bf(A3)} hold,  $\psi^*(x)=\psi(x)=c_1x^2+c_2|x|^{\alpha}$,
	and $\psi\in\mathrm{WLSC}(\alpha,1)$. According to \cite[(12)]{GR17}, we have $$H(x)\leqslant \frac{10}{\pi(\alpha-1)}\frac{1}{c_1{|x|}^{-1}+c_2|x|^{1-\alpha}}.$$
	Combining it with Corollary \ref{wlsc},   if 
	$$\sup_{x}|x|^{\alpha-1} \mu((-|x|,|x|)^c)<\infty,$$
	then  ${Y}$ is exponentially ergodic, and the spectral gap  $$\lambda_{1}\geqslant\frac{\pi(\alpha-1)}{80}\inf_x \frac{c_1{|x|}^{-1}+c_2|x|^{1-\alpha}}{\mu((-|x|,|x|)^c)};$$
	if 		\begin{equation*}
	\int_{\mathbb{R}}\sigma(x)^{-\alpha}|x|^{\alpha-1} d x<\infty,
	\end{equation*}
	then   ${Y}$ is strongly  ergodic, and the convergence rate in strong ergodicity
	$$\kappa\geqslant\frac{\pi(\alpha-1)}{20\int_{\mathbb{R}}(a(x))^{-1}\left(c_1{|x|}^{-1}+c_2|x|^{1-\alpha}\right)^{-1}\d x}.$$
	\end{example}
\begin{remark}
	Note that if $c_1=0$ (resp., $c_2=0$), then the result is reduced to the time-changed symmetric stable process (resp., time-changed Brownian motion). 
\end{remark}
\begin{example}
	Let $\psi(x)=x^2+|x|$.  The process $X$ associated with $\Delta+\Delta^{1/2}$ is a sum of the Cauchy process and independent Brownian motion. Obviously, {\bf(A1)}--{\bf(A3)} hold. By \cite[Lemma 2.14]{GR17},  
	$$H(x)\leqslant\frac{10}{\pi}\int_{1/|x|}^\infty\frac{\d r}{r+r^2}=\frac{10}{\pi}\log\left(1+|x|\right).$$
	Therefore, if
	$$\delta_2:=\sup_{x} \log(1+|x|)\mu((-|x|,|x|)^c)<\infty,$$
	then $Y$ is exponentially ergodic, and the spectral gap $$\lambda_{1}\geqslant\frac{\pi}{80\delta_2}.$$
	
	If	$$I_2:=\int_{\mathbb{R}}a(x)^{-1}\log(1+|x|)d x<\infty,$$ then ${Y}$ is strongly  ergodic and $$\kappa\geqslant\frac{\pi}{20I_2}.$$
	
\end{example}

\section{Time change and Green potential }\label{prem}

Let $X=\left(X_{t}\right)_{t \geqslant 0}$ be a one-dimensional  symmetric  L\'{e}vy process with L\'{e}vy measure $\nu$, 
transition density $p_t(x,y)=p_t(y-x)$, and characteristic exponent $\psi$ given by \eqref{symbol}.

Recalling that $a$ is a positive and locally bounded measurable function on $\mathbb{R}$ with $a(x)^{-1}$ is Lebesgue integrable, 
$A_{t}=\int_{0}^{t} 1 / a(X_{s}) d s$ is the positive continuous additive functional and $Y$ is the time-changed L\'{e}vy process defined as \eqref{time change}. The Revuz measure $\mu$ of $A_t$ with respect to $\d x$, is given by (cf. \cite{FG88})
$$\mu(f)=\lim_{t\rightarrow\infty}\frac{1}{t}\int_\mathbb{R}\mathbb{E}_x\left[\int_0^tf(X_s)\d A_s\right]\d x.$$
Since $\d x$ is the invariant measure of $X$,  for nonnegative bounded function $f$, we have
\begin{equation}
	\begin{split}
		\mu(f)&=\lim_{t\rightarrow\infty}\frac{1}{t}\int_\mathbb{R}\mathbb{E}_x\left[\int_0^tf(X_s)a(X_s)^{-1}\d s\right]\d x=\lim_{t\rightarrow\infty}\frac{1}{t}\int_0^t\int_\mathbb{R}P_s(fa^{-1})(x)\d x\d s\\
		&=\lim_{t\rightarrow\infty}\frac{1}{t}\int_0^t\int_\mathbb{R}(fa^{-1})(x)\d x\d s=\int_\mathbb{R}(fa^{-1})(x)\d x,
	\end{split}
\end{equation}
thus the Revuz measure $\mu(\d x)=a(x)^{-1}\d x$. 
Combining this fact and \cite[Theorem 5.2.2, Theorem 5.2.8 and Corollary 5.2.12]{CM12}, similar to \cite[Page 2807]{CW14}, we know that $Y$ is $\mu$-symmetric and  
its Dirichlet form $(\widehat{\scr{E}},\widehat{\mathcal{F}})$ is given by 
\begin{equation}\label{tc-Diri}
	\widehat{\scr{E}}(f,g)=\scr{E}(f,g),\quad f,g\in \widehat{\mathcal{F}}:=\mathcal{F}_e\cap L^2(\mu),
\end{equation}
where $\mathcal{F}_e$ is the extended Dirichlet space of $(\scr{E},\mathcal{F})$, i.e., the family of  functions $u$ satisfy that there exists an $\scr{E}$-Cauchy sequence $\{u_n\}\subset \mathcal{F}$ such that for a.e. $x$, $\lim_{n \rightarrow \infty}u_n=u$ in $L^2(\d x)$ and $\scr{E}(u,u)=\lim_{n\rightarrow\infty}\scr{E}(u_n,u_n)$.
Therefore, the $L^{2}$ infinitesimal generator of $Y$ is $\mathcal{L}=a \mathcal{A}$. 

By a similar argument to \cite[Corollary 4.2]{bsw13}, we can also prove that the extended infinitesimal generator $\widetilde{\mathcal{L}}$ (see \cite[Definition 2.1]{st05}),  is also $a \mathcal{A}$.




Recalling that $P_t^0$ is the killed semigroup of $X$ given by \eqref{killed semigroup}. Define the \textit{Green potential measure} $G_X^0(x,A)$ for $P_t^{0}$ by
$$G_X^0(x,A)=\int_{0}^\infty P_t^0(x,A)\d t,\quad \forall x\in \mathbb{R}, \ \forall A\in\mathscr{B}(\mathbb{R}).$$
Denote by $P_t^{0,Y}$ the  semigroup of $Y$ killed upon hitting the origin, i,e. $$P_t^{0,Y}(x,A)=\mathbb{P}_x[Y_t\in A,t<\tau_0^Y],$$ where $\tau_0^Y=\inf\{t>0:Y_t\neq0\}$.
Let $G_Y^0(x,A)$ be the Green measure of $Y$ killed upon $0$:
$$G_Y^0(x,A)=\int_{0}^\infty P_t^{0,Y}(x,A)\d t,\quad \forall x\in \mathbb{R}, \ \forall A\in\mathscr{B}(\mathbb{R}).$$
Similar to \cite[(14)]{W21}, for the time-changed process $Y$, 
\begin{equation}\label{tc property}
	G_Y^0f(x)= \int_{\mathbb{R}}G_X^0(x,y)a(y)^{-1}\d y.\\
\end{equation}
\section{Proofs of main results}
To finish the proof of Theorem \ref{main thm}, we need to consider
 the first Dirichlet eigenvalue
\begin{equation}\label{vari diri}
	\lambda_{0}=\inf \{\mathscr{E}(f, f): f \in \widehat{\mathcal{F}}, \mu(f^{2})=1 \text { and }f(0)=0\},
\end{equation}
which will play a crucial role in the proof of Theorem \ref{main thm}.

First, we introduce the following dual variational inequality, which is mainly motivated by the dual variational formulas for one-dimensional diffusion processes (see \cite[Theorem 6.1]{cmf05}).

\begin{lemma}\label{variational}
	Let $G_Y^{0}$ be the Green operator of $Y$ killed upon $\{0\}$, $C(B)$ be the space of all continuous functions on a measurable set $B\subset\mathbb{R}$. Denote by
	$$
	\begin{aligned}
		&\mathscr{H}=\left\{f: f(0)=0, f \in C( \mathbb{R}\setminus\{0\})\right\} \\
		&\widetilde{\mathscr{H}}=\left\{f:  f(0)=0, \text { there exists } x_{0} >0 \right. \text { such that }\left.f=f\left(\cdot \wedge x_{0}\vee (-x_0)\right), f \in C\left(-x_0, x_{0}\right)\right\}.
	\end{aligned}
	$$
	Then $$\inf_{f\in \widetilde{\mathscr{H}}} 
	\sup_{x \in \mathbb{R}\setminus\{0\}} \frac{f(x)}{G_Y^{0} f(x)}\geqslant\lambda_{0}\geqslant\sup_{f\in \mathscr{H}} 
	\inf_{x \in \mathbb{R}\setminus\{0\}} \frac{f(x)}{G_Y^{0} f(x)}.$$
\end{lemma}

\prf
First we proof the upper bound. Note that for any $g\in \widetilde{\mathscr{H}}$, there exists $x_0> 0$ such that $g(-x_0)\leqslant g\leqslant g( x_{0})$. Therefore, $g\in L^2(\pi)$. By a similar argument to the proof of \cite[Theorem 2]{W21} (see \cite[Page 12]{W21}),  we have	 
$G_Y^{0}g\in\widehat{\mathcal{F}}$, and 
$$\scr{E}(G_Y^{0}g,G_Y^{0}g)=\int gG_Y^{0}g \d \pi.$$ 
By the definition \eqref{vari diri},
$$\lambda_{0}\leqslant \frac{\scr{E}(G_Y^{0}g,G_Y^{0}g)}{\pi((G_Y^{0}g)^2)}\leqslant \sup_{x \in \mathbb{R}\setminus\{0\}} \frac{g(x)}{G_Y^{0} g(x)}.$$
Now we get the upper bound by the arbitrariness of $g\in\widetilde{\mathscr{H}}.$

Next, denote by $$\lambda_{0}^{(n)}=\inf\{\scr{E}(f,f):\mu(f^2)=1,f|_{(-\infty,-n)\cup(n,\infty)}=0\}.$$
According to \cite[Lemma 11]{W21},
$$\lim_{n \rightarrow \infty}\lambda_{0}^{(n)}=\lambda_{0}.$$
By a similar argument to the proof of \cite[Lemma 7]{W21}, for regular set $B_n:=(-n,0)\cup(0,n)$  (see \cite[Page 68]{chung86} for the definition of regular set), 
$$\lambda_{0}^{(n)}  \geqslant \sup_{f\in C_b(B_n)}
\inf _{x \in B_n} \frac{f(x)}{G_Y^{B_n} f(x)},$$
where $G_Y^{B_n}$ is the Green operator defined as
$G_Y^{B_n}(x,A):=\int_{0}^\infty\mathbb{P}_x[Y_t\in A,t<\tau_{B_n}^Y]\d t,$
$\tau_{B_n}^Y$ is the exit time from $B_n$: $\tau_{B_n}^Y:=\inf\{t\geqslant0: Y_t\notin B_n\}$.

Next, since $G_Y^{B_n} \varphi\leqslant G_Y^{\{0\}^c} \varphi$ and $f\in C(\mathbb{R}\setminus\{0\})$ is bounded on $B_n$,   we have that for any $f\in \mathscr{H}$,
\begin{equation}
	\begin{split}
		\lambda_{0}&=\lim _{n\rightarrow \infty}\lambda_{0}^{(n)}  \geqslant\lim _{n\rightarrow \infty}
		\inf _{x \in B_n} \frac{f(x)}{G_Y^{0} f(x)}\geqslant 
		\inf _{x \in \mathbb{R}\setminus\{0\}} \frac{f(x)}{G_Y^{0} f(x)},
	\end{split}
\end{equation}
thus we obtain the lower bound.
\deprf

The following explicit estimates for lower and upper bounds of $\lambda_0$ is similar to \cite[Theorem 2]{W21}, and can be obtained directly by Lemma \ref{variational}.
\begin{theorem}
	Let $\delta$ be defined by \eqref{delta} and 
	$$\delta^+=\sup_{x>0} H(x)\mu((x,\infty)), \quad \delta^-=\sup_{x<0} H(x)\mu((-\infty,-x)).$$ Then
	\begin{equation}
		\frac{1}{\delta_+}+\frac{1}{\delta_-}\geqslant\lambda_{0}\geqslant\frac{1}{8\delta}.
	\end{equation}
\end{theorem}
\prf First, we prove the upper bound. 
According to Lemma \ref{variational}, for any  $g\in\widetilde{\mathscr{H}}$,
\begin{equation}\label{upbd}
	\lambda_{0}\leqslant\sup_{x \in \mathbb{R}\setminus\{0\}} \frac{g(x)}{G_Y^{0} g(x)}\leqslant \sup_{x>0}\frac{g(x)}{G_Y^{0}g(x)}+\sup_{x<0}\frac{g(x)}{G_Y^{0}g(x)}.
\end{equation}
Note that by \cite[Proposition 2.4]{GR17}, for $xy>0$, $$G_Y^{0}(x,y)\geqslant H(|x|\wedge|y|), $$
so we get that for $x>0$,
\begin{equation}
	{G_Y^{0}g(x)}\geqslant \int_{0}^{\infty}H(x\wedge y)g(y)a(y)^{-1} \d y\geqslant H(x)\int_{x}^\infty g(y)a(y)^{-1} \d y.
\end{equation}
Since $g\in\widetilde{\mathscr{H}}$, there exists $x_0>0$ such that
$g(x)=	g(x\wedge x_0\vee (-x_0))$. Now by choosing $g(x)=H(x\wedge x_0\vee (-x_0))$, we obtain that for $x>0$,
\begin{equation}\label{G+}
	\frac{G_Y^{0}g(x)}{g(x)}=\frac{G_Y^{0}g(x\wedge x_0)}{g(x\wedge x_0)}\geqslant H(x_0)\int_{x_0}^\infty a(y)^{-1} \d y.
\end{equation}
Similarly, for $x<0$, 
\begin{equation}\label{G-}
	\frac{G_Y^{0}g(x)}{g(x)}=\frac{G_Y^{0}g(x\vee (-x_0))}{g(x\vee (-x_0))}\geqslant H(-x_0)\int_{-\infty}^{-x_0}a(y)^{-1} \d y.
\end{equation}
Therefore, by combining  \eqref{upbd}, \eqref{G+}, and \eqref{G-}, we obtain the estimate for upper bound.

Next, we consider the lower bound. Let $G^0$ be the Green operator of $Y$ killed upon $\mathbb{R}\setminus\{0\}$. According to \cite[Proposition 2.3 and 2.4]{GR17}, 
\begin{equation}\label{G0-general}
	G_X^0(x, y)=H(x)+H(y)-H(y-x)\leqslant 2 (H(x)\wedge H(y)).
\end{equation} 
By the property of time change and \eqref{G0-general}, for any $f$ with $	G_Y^{0} |f|<\infty$,  
\begin{equation}\label{II-general}
	G_Y^0f(x)
	\leqslant \int_{\mathbb{R}}2 (H(y)\wedge H(x))f(y)a(y)^{-1}\d y.\\
\end{equation}
Since 	$\psi(t)/t\rightarrow\infty$ as $t\rightarrow\infty$, by Dirichlet criterion,  
$ \int_{0}^{\infty} t\sin (x t) \psi(t)^{-1}  dt$
is integrable. Therefore, $H$ is differentiable, and
$$H'(x)=\frac{1}{\pi} \int_{0}^{\infty} \frac{t\sin x t}{\psi(t)} d t<\infty.$$
Recalling that  $p_t(x)=(2\pi)^{-1}\int_{\mathbb{R}}\e^{-t\psi(\xi)-i\xi x}\d \xi$ and
$\psi(\xi)=\psi(-\xi)$. Since $$\int_{\mathbb{R}}\e^{-t\psi(\xi)}\d \xi=p_t(0)=p_t(x,x)<\infty,$$
$\e^{-t\psi}\in L^1(\d x)$.
By \cite[Theorem 1.1]{GT13}, for $|\eta|=r$,
$$-\frac{1}{2\pi r}\frac{\d}{\d r}p_t(r)=\int_{\mathbb{R}^3}\e^{-t\psi(|x|)}\e^{-2\pi ix\cdot\eta}\d x=\frac{1}{2\pi}q_t(\eta),$$
where $q_t(\eta)$ is the transition density of a 3-dimensional L\'{e}vy process with characteristic exponent $\Psi(\eta)=\psi(|\eta|)$. 
Therefore, $q_t(\eta)\geqslant 0$,  $p_t(x)$ is non-increasing for $x$, i.e. it is unimodal (see \cite{GR17}). Thus by \cite{GR17}, $H(x)$ is non-increasing on $(0,\infty)$ and non-decreasing on $(-\infty,0)$, so we have
\begin{equation}
	\begin{split}		
		G_Y^0f(x)
		&\leqslant2 \left(\int_{\mathbb{R}\setminus(-|x|,|x|)} H(x)f(y)\mu(d y)+\int_{-|x|}^{|x|} H(y)f(y)\mu(d y)\right)\\
		&=2\int_{0}^{|x|} H'(z)\left(\int_{\mathbb{R}\setminus(-z,z)} f(y)\mu(d y)\right)d z.\\
	\end{split}
\end{equation}
Choose $f(x)=\sqrt{H(x)}$.
By using integration by parts, for any $y>0$,
\begin{equation*}
	\begin{split}
		\int_{\mathbb{R}\setminus(-y,y)} \sqrt{H(z)}\mu(d z)
		\leqslant 	&\sqrt{H(y)}\mu((-y,y)^c)\\
		&+\int_{y}^{\infty}\frac{H'(z)\mu((-z,z)^c)}{2\sqrt{H(z)}}d z.\\
	\end{split}
\end{equation*}
Since $\delta<\infty$, then 
\begin{equation*}
	\begin{split}
		\int_{\mathbb{R}\setminus(-y,y)} \sqrt{H(z)}\mu(d z)
		\leqslant\frac{\delta}{\sqrt{H(y)}}+\frac{\delta}{2}\int_{y}^{\infty}{H(z)}^{-3/2}d H(z)
		=\frac{2\delta}{\sqrt{H(y)}}.\\
	\end{split}
\end{equation*}
Thus 
\begin{equation*}
	\frac{G_Y^0\sqrt{H}}{\sqrt{H}}(x)\leqslant\frac{2}{\sqrt{H(x)}}\int_{0}^{|x|} H'(z)\frac{2\delta}{\sqrt{H(z)}}d z={8\delta}.
\end{equation*}	

Now by Lemma \ref{variational} and letting $f=\sqrt{H}$, we obtain that if $\delta<\infty$, then
$$\lambda_{0}\geqslant  \inf_{x\neq0} \frac{\sqrt{H(x)}}{G_Y^{0}\sqrt{H(x)}}     \geqslant\frac{1}{8\delta}.$$
\deprf


\noindent{\bf Proof of Theorem \ref{main thm}}. 

(1)
By \cite[Proposition 3.2]{cmf00'}, $\lambda_{1}\geqslant\lambda_{0}$. Hence  
$$\lambda_1\geqslant\lambda_{0}\geqslant 
\frac{1}{8\delta},$$ and ${Y}$ is exponential ergodicity.

(2) Specially, by choosing $f\equiv 1$ in \eqref{II-general}, we have
$$M_0:=\sup_x\mathbb{E}_x\tau_{0}^Y=\sup_xG_Y^01(x)
\leqslant \int_{\mathbb{R}}2 H(y)a(y)^{-1}d y<\infty,\\$$
thus by \cite[Theorem 1.2(R2)]{mwt20},
$\kappa\geqslant M_0^{-1}\geqslant(2I)^{-1}>0.$\deprf

\noindent{\bf Proof of Corollary \ref{wlsc}}.

First, according to \cite[Lemma 2.14]{GR17}, $H(x)\approx(|x|\psi(1/x))^{-1}$, thus by using  Theorem \ref{main thm}, we obtain the exponential ergodicity and strong ergodicity.

By \cite[Lemma 2.14 and (12)]{GR17},
$$H(x)\leqslant\frac{10}{\pi\beta}\int_{1/x}^\infty\frac{1}{\psi^*(s)}d s\leqslant \frac{10}{\pi\beta^2(\alpha-1)x\psi^*(1/x)},\ \text{for}\ x>0.$$
Then the estimates for $\lambda_{1}$ and $\kappa$ follow from Theorem \ref{main thm}. 

\deprf
\section*{Acknowledgements}
The author thanks Prof.\ Yong-Hua Mao and Prof.\ Jian Wang for valuable conversations of this paper. The author would also like to thank Mr.\ Zhi-Feng Wei for his helpful suggestions. This work is supported in part by the  National Key Research and Development Program of China (2020YFA0712900), the National Natural Science Foundation of China (Grant No.11771047), and the project from the Ministry of Education in China.

\bibliographystyle{plain}
\bibliography{symlevy}

\begin{thebibliography}{10}

\bibitem{bsw13}
B.~Bottcher, R.~Schilling, and J.~Wang.
\newblock {\em {L\'{e}vy-type processes: construction, approximation and sample
  path properties. From: L\'{e}vy matters III, volume 2099 of Lecture Notes in
  Mathematics, with a short biography of Paul L\'{e}vy by Jean Jacod.}}
\newblock Berlin: Springer, 2013.

\bibitem{cmf00'}
M.~F. Chen.
\newblock {Explicit bounds for the first eigenvalues.}
\newblock {\em Sci. in China Ser. A}, 43:1051--1059, 2000.

\bibitem{cmf05}
M.F. Chen.
\newblock {\em {Eigenvalues, inequalities, and ergodic theory}}.
\newblock London: Springer, 2004.

\bibitem{CW14}
Z.~Q. Chen and J.~Wang.
\newblock Ergodicity for time-changed symmetric stable processes.
\newblock {\em Stoch. Proc. Appl.}, 124(9):2799--2823, 2014.

\bibitem{CM12}
Z.Q. Chen and M.~Fukushima.
\newblock {\em {Symmetric Markov processes, time change, and boundary theory}}.
\newblock Princeton: Princeton Univ. Press, 2012.

\bibitem{chung86}
K.L. Chung.
\newblock {Doubly-Feller process with multiplicative functional}.
\newblock {\em In Seminar on Stochastic Processes, Birkh\"{a}user.}, 12:63--78,
  1986.

\bibitem{FG88}
P.~J. Fitzsimmons and R.~K.Getoor.
\newblock {Revuz measures and time changes}.
\newblock {\em Math. Z.}, 199(2):233--256, 1988.

\bibitem{GT13}
L.~Grafakos and G.~Teschl.
\newblock {On Fourier transforms of radial functions and distributions}.
\newblock {\em J. Fourier Anal. Appl}, 19(1):167--179, 2013.

\bibitem{GR17}
T.~Grzywny and M.Ryznar.
\newblock {Hitting times of points and intervals for symmetric L\'{e}vy
  processes}.
\newblock {\em Potential Anal.}, 46:739--777, 2017.

\bibitem{wj19}
D.~Luo and J.~Wang.
\newblock {Refined basic couplings and Wasserstein-type distances for SDEs with
  Levy noises}.
\newblock {\em Stoch. Proc. Appl.}, 129:3129–--3173, 2019.

\bibitem{myh06}
Y.H. Mao.
\newblock {Convergence rates in strong ergodicity for Markov processes}.
\newblock {\em Stoch. Proc. Appl.}, 116:1964--1976, 2006.

\bibitem{mwt20}
Y.H. Mao and T.~Wang.
\newblock {Convergence rates in strong ergodicity by hitting times and
  $L^2$-exponential convergence rates}.
\newblock {\em see ArXiv 2102.07069}.

\bibitem{OY13}
Y.~Oshima.
\newblock {\em {Semi-Dirichlet forms and Markov processes}}.
\newblock Berlin: De Gruyter Studies in Mathematics., 2013.

\bibitem{DZ96}
G.~Da Prato and J.~Zabczyk.
\newblock {\em {Ergodicity for infinite-dimensional systems}}.
\newblock Cambridge: Cambridge University Press, 1996.

\bibitem{st05}
Y.~Shiozawa and M.~Takeda.
\newblock {Variational formula for Dirichlet forms and estimates of principal
  eigenvalues for symmetric $\alpha$-stable processes}.
\newblock {\em Potent. Anal.}, 23:135--151, 2005.

\bibitem{wj08}
J.~Wang.
\newblock Criteria for ergodicity of l\'{e}vy type operators in dimension one.
\newblock {\em Stoch. Proc. Appl.}, 118:1909--1928, 2008.

\bibitem{wj13}
J.~Wang.
\newblock {Exponential ergodicity and strong ergodicity for SDEs driven by
  symmetric $\alpha$-stable processes.}
\newblock {\em Appl. Math. Lett.}, 26:654--658, 2013.

\bibitem{W21}
T.~Wang.
\newblock Exponential and strong ergodicity for one-dimensional symmetric
  stable jump diffusions.
\newblock {\em See ArXiv: 2104.07947}, 2021.

\bibitem{ya10}
K.~Yano.
\newblock {Excursions away from a regular point for one-dimensional symmetric
  L\'{e}vy processes without Gaussian part}.
\newblock {\em Potential Anal.}, 32(4):305--341, 2010.

\end{thebibliography}

\end{document}